\documentclass[12pt]{article}
\usepackage{amssymb,amsmath}
\usepackage[amsmath,thmmarks,thref]{ntheorem}
\usepackage{graphics}
\usepackage{enumerate}
\usepackage{dsfont}
\usepackage{mathrsfs}
\makeatletter
\def\s@btitle{\relax}
\def\subtitle#1{\gdef\s@btitle{#1}}
\def\@maketitle{%
  \newpage
  \null
  \vskip 2em%
  \begin{center}%
  \let \footnote \thanks
    {\LARGE \@title \par}%
                \if\s@btitle\relax
                \else\typeout{[subtitle]}%
                        \vskip .5pc
                        \begin{large}%
                                \textsl{\s@btitle}%
                                \par
                        \end{large}%
                \fi
    \vskip 1.5em%
    {\large
      \lineskip .5em%
      \begin{tabular}[t]{c}%
        \@author
      \end{tabular}\par}%
    \vskip 1em%
    {\large \@date}%
  \end{center}%
  \par
  \vskip 1.5em}
\makeatother

\input{diagxy}
\xyoption{curve}
\newbox\anglebox 
\setbox\anglebox=\hbox{\xy \POS(75,0)\ar@{-} (0,0) \ar@{-} (75,75)\endxy}
 \def\pbangle{\copy\anglebox}
\newbox\angleboxr 
\setbox\angleboxr=\hbox{\xy \POS(0,0)\ar@{-} (0,75) \ar@{-} (75,0)\endxy}
 \def\pbangler{\copy\angleboxr}
\newbox\sanglebox 
\setbox\sanglebox=\hbox{\xy \POS(50,0)\ar@{-} (0,0) \ar@{-} (50,50)\endxy}
 
\newbox\sangleboxr 
\setbox\sangleboxr=\hbox{\xy \POS(0,0)\ar@{-} (0,50) \ar@{-} (50,0)\endxy}
 
\newbox\sangleboxf 
\setbox\sangleboxf=\hbox{\xy \POS(50,50)\ar@{-} (50,0) \ar@{-} (0,50)\endxy}
 
\newbox\angleboxf 
\setbox\angleboxf=\hbox{\xy \POS(75,75)\ar@{-} (75,0) \ar@{-} (0,75)\endxy}
 
\newbox\sangleboxfr 
\setbox\sangleboxfr=\hbox{\xy \POS(0,50)\ar@{-} (50,50) \ar@{-} (0,0)\endxy}
 
\newbox\angleboxfr 
\setbox\angleboxfr=\hbox{\xy \POS(0,75)\ar@{-} (75,75) \ar@{-} (0,0)\endxy}

\newdir{|>}{!/4.7pt/\dir{|}
        *:(1,-.2)\dir^{>}
        *:(1,+.2)\dir_{>}}

\def\embedd{\to/^{ (}->/}
\newcommand{\pair}[1]{\ensuremath{\langle {#1} \rangle}}
\newcommand{\sublat}[2]{\ensuremath{\operatorname{Sub}_{#1}\left({#2}\right)}}
\newcommand{\cat}[1]{\ensuremath{\mathcal{#1}}}
\newcommand{\topo}[1]{\ensuremath{\mathscr{#1}}}
\newcommand{\thry}[1]{\ensuremath{\mathbb{#1}}}
\newcommand{\synt}[2]{\ensuremath{\mathcal{#1}_{\mathbb{#2}}}}
\newcommand{\alg}[1]{\ensuremath{\mathbf{#1}}}
\newcommand{\mng}[1]{\ensuremath{\mathrm{#1}}}

\newcommand{\Sets}{\ensuremath{\mathbf{Set}}}

\newcommand{\Sh}[1]{\protect\ensuremath{\operatorname{Sh}\left(#1\right)}}

\newcommand{\theory}{\ensuremath{\mathbb{T}}}

\newcommand{\cterm}[2]{\ensuremath{\left \{ {#1}\ \; \vrule \; \ {#2}\right \}}}
\newcommand{\fratil}[3]{\ensuremath{#1 \! :\! #2\to<125> #3}}
\newcommand{\inklusjon}[3]{\ensuremath{#1 \! :\! #2\to/^{ (}->/<150> #3}}

\newcommand{\classtop}[1]{\ensuremath{\mathbf{Set}[\thry{#1}]}}

\newcommand{\Eqsheav}[2]{\protect\ensuremath{\operatorname{Sh}_{#1}(#2)}}
\begin{document}
%
\title{Subgroupoids and Quotient Theories}
\author{Henrik Forssell
\thanks{Thanks to Ji\v{r}\'{i} Rosick\'{y} and Steve Awodey for interesting and helpful conversations. This research was supported by the Eduard \v{C}ech Center for Algebra and Geometry, grant no.\ LC505.}}
\maketitle
\begin{abstract}
Moerdijk's site description for equivariant sheaf toposes on open topological groupoids is used to give a proof for the (known, but apparently unpublished) proposition that
if \topo{H} is a subgroupoid of an open topological groupoid \topo{G}, then  the topos of equivariant sheaves on \topo{H}  is a subtopos of the topos of equivariant sheaves on \topo{G}. This proposition is then applied to the study of quotient geometric theories and subtoposes. In particular, an intrinsic characterization is given of those subgroupoids that are definable by quotient theories.
\end{abstract}
%
%
%
%
\section{Introduction}
\label{Section> Introduction}
In \cite{butz:98b}, Butz and Moerdijk showed that a topos with enough points can be represented as the topos of equivariant sheaves on an open topological groupoid constructed from points of the topos. In `logical' terms, this can be rephrased as saying that for any geometric theory \theory\ with enough models, there exists an open topological groupoid \topo{G} consisting of \theory-models and isomorphisms such that the classifying topos of \theory\ is equivalent
to the topos of equivariant sheaves on \topo{G}
\begin{equation}\label{Equation: Introequivalence}\classtop{T} \simeq\Eqsheav{G_1}{G_0}\end{equation}
Conversely, any equivariant sheaf topos \Eqsheav{G_1}{G_0} classifies a geometric theory with enough models, and \topo{G} can be regarded as consisting of \theory-models and isomorphisms.
Considering the displayed equivalence (\ref{Equation: Introequivalence}), there is on the `theory' side a correspondence between subtoposes of \classtop{T} and quotient theories of \theory\ (see \cite[Theorem 3.6]{caramello:99LoT}). On the groupoid side, it is known to specialists (Moerdijk in particular) that a subgroupoid of an open topological groupoid induces a subtopos of equivariant sheaves, but this fact appears not to have been published. As a first outline of the connection between subgroupoids of \topo{G} and quotient theories of \theory, this paper  first fills in a proof of that fact and points out the resulting Galois connection between subgroupoids of \topo{G} and subtoposes of \Eqsheav{G_1}{G_0} (and thus quotient theories of \theory), and then characterizes the subgroupoids of \topo{G} that are definable by quotient theories. The whole investigation is carried out using Moerdijk's site description for equivariant sheaf toposes given in \cite{moerdijk:88}, and a brief introduction to that construction is given first.

\section{Subgroupoids and subtoposes}
\label{Section: Subgroupoids and subtoposes}
\subsection{Groupoids and equivariant sheaves}
\label{Subsection: Groupoids and equivariant sheaves}
This section briefly recalls from \cite{elephant1}, \cite{moerdijk:88}, \cite{moerdijk:90} the topos of equivariant sheaves on a topological groupoid and Moerdijk's site description for such toposes (written out here for topological rather than localic groupoids and writing out a few additional details, cf.\ especially \cite[\S 6]{moerdijk:88}, a more detailed and self-contained presentation can be found in the online note \cite{onlinenote:subgpds}). Let \topo{G} be a topological groupoid, fully written out as a groupoid object in the category \alg{Sp} of topological spaces and continuous maps as
\[\bfig
%
%
\morphism|m|<750,0>[G_1\times_{G_0} G_1`G_1;m]
%
%
\morphism(750,0)|a|/@{>}@<5pt>/<750,0>[G_1`G_0;d]
\morphism(750,0)|m|/@{<-}/<750,0>[G_1`G_0;e]
\morphism(750,0)|b|/@{>}@<-5pt>/<750,0>[G_1`G_0;c]
\Loop(750,0)G_1(ur,ul)_{\mng{i}}
%
\efig\]
with $m$ the composition, $e$ the mapping to identities, and $i$ the mapping to inverses. This notation will be mixed with the usual notation $g\circ f$, $1_x$, $f^{-1}$. \topo{G} is called \emph{open} if the domain and codomain maps are open. It follows that composition of arrows must also be open. The objects of the category of \emph{equivariant sheaves}, \Eqsheav{G_1}{G_0}, on \topo{G} are pairs \pair{r:R\rightarrow G_0,\rho} where $r$ is a local homeomorphism---i.e.\ an object of \Sh{G_0}---and $\rho$ is a continuous action, i.e.\ a continuous map
\[\rho:G_{1}\times_{G_{0}}R\to R\]
with the pullback being along the domain map and such that $r(\rho(f,x))=c(f)$, satisfying the expected unit and composition axioms.
%
If \topo{G} is an open topological groupoid then it follows that the action $\rho$ is an open map.
A morphism of equivariant sheaves is a morphism of sheaves (local homeomorphisms) commuting with the actions. The category, \Eqsheav{G_1}{G_0}, of equivariant sheaves on \topo{G} is a (Grothendieck) topos. The forgetful functors of forgetting the action, $\fratil{u}{\Eqsheav{G_1}{G_0}}{ \Sh{G_0}}$, and of forgetting the topology, $\fratil{v}{\Eqsheav{G_1}{G_0}}{ \Sets^{\topo{G}}}$, are both conservative inverse image functors. A \emph{continuous functor}, or morphism of topological groupoids, $\fratil{f}{\topo{H}}{\topo{G}}$, i.e.\ a morphism of groupoid objects in \alg{Sp}
%
%
%
%
%
%
%
induces a geometric morphism $\fratil{f}{\Eqsheav{H_1}{H_0}}{\Eqsheav{G_1}{G_0}}$
%
%
%
%
%
%
%
%
where $f^*$ pulls a sheaf back along $f_0$ and equips it with an action using $f_1$ in the expected way (both $u^*$ and $v^*$ above are examples).
%
%
%
%
%
%
%
%

Let \topo{G} be an open topological groupoid. Let $N\subseteq G_1$ be an open subset closed under composition and inverse, and let $U=d(N)=c(N)$. Refer to such a pair $(U,N)$ as an \emph{open subgroupoid}. Form the quotient space $d^{-1}(U)\twoheadrightarrow d^{-1}(U)/_{\sim_N}$ by $f\sim_N g$ iff $c(f)=c(g)$ and $g^{-1}\circ f \in N$. The quotient map $q$ is then an open surjection, the codomain map $c:d^{-1}(U)/_{\sim_N}\rightarrow G_0$ is a local homeomorphism, and composition defines a continuous action on $d^{-1}(U)/_{\sim_N}$, so that we have an equivariant sheaf denoted \pair{\topo{G},U,N}.
Objects of the form $\pair{\topo{G},U,N}$ form a generating set for \Eqsheav{G_1}{G_0}. Briefly, given an equivariant sheaf \pair{r:R\rightarrow G_0,\rho} and a continuous section $t:U\rightarrow R$, we get an open  set of arrows
\[N_t=\cterm{f\in d^{-1}(U)\cap c^{-1}(U)}{\rho(f,t(d(f)))=t(c(f))}\]
(by pullback of the open set $t(U)$ along an appropriate continuous map) which is closed under composition and inverse, and such that $d(N)=c(N)=U$.
%
%
There is a canonical continuous section $e:U\rightarrow d^{-1}(U)/_{\sim_N}$ defined by $x\mapsto [1_x]_{\sim_{N_t}}$, and the section $t$ lifts to a morphism, $\hat{t}:\pair{\topo{G},U,N_t}\rightarrow R$, of \Eqsheav{G_1}{G_0},
\begin{equation}\label{Equation: Nt and t-hat}\bfig
\Vtriangle/>`<-`<-/<350,250>[d^{-1}(U)/_{\sim_{N_t}}`R`U;\hat{t}`e`t]
\place(350,150)[=]
\efig\end{equation}
such that $\hat{t}([f])=\rho(f,t(d(f)))$. One easily sees that $\hat{t}$ is 1-1. For reference:
%
\begin{proposition}\label{Proposition: GUNS generate}
Any object $A\in\Eqsheav{G_1}{G_0}$ is the join of its subobjects of the form $\pair{\topo{G},U,N}\rightarrowtail A$ for open subgroupoids $(U,N)$.
%
%
%
%
\end{proposition}

The full subcategory of \Eqsheav{G_1}{G_0} consisting of objects of the form \pair{\topo{G},U,N} is, accordingly, a site for \Eqsheav{G_1}{G_0} when equipped with the canonical coverage. Refer to this as the \emph{Moerdijk site} for \Eqsheav{G_1}{G_0}, and denote it $\cat{S}_{\topo{G}}\to/^{ (}->/<125> \Eqsheav{G_1}{G_0}$. Moerdijk sites are closed under subobjects. For consider an object   \pair{\topo{G},U,N} and let $V\subseteq U$ be an open subset closed under $N$, that is, such that $x\in V$ and $f:x\rightarrow y$ in $N$ implies $y\in V$. Then
\[d^{-1}(V)/_{\sim_{N\upharpoonright_{V}}}= m(G_1\times_{G_0} e(V))\subseteq d^{-1}(U)/_{\sim_{N}}\]
is an open subset closed under the action, and so a subobject. All subobjects are of this form:
\begin{lemma}\label{Lemma: Subobjects in the Moerdijk site}
Let \pair{\topo{G},U,N} be an object of $\cat{S}_{\topo{G}}$.  Then
$V \mapsto d^{-1}(V)/_{\sim_{N\upharpoonright_{V}}}$
defines an isomorphism between the frame of open subsets of $U$ that are closed under $N$ and the frame of subobjects of  \pair{\topo{G},U,N}.
\begin{proof}The inverse is given by pulling back along the canonical section $e: U\rightarrow d^{-1}(U)/_{\sim_{N}}$.
\end{proof}
\end{lemma}
The morphisms in the Moerdijk site can be described in a manner similar to the objects in it. Consider a morphism $\hat{t}:\pair{\topo{G},U,N}\rightarrow \pair{\topo{G},V,M}$.
It is easily seen that such a morphism determines and is determined by a section $t:U\rightarrow d^{-1}(V)/_{\sim_M}$ with the property that for any $f:x\rightarrow y$ in $N$, we have that $f\circ t(x)=t(y)$. And such a section can be described as an open set:
\begin{lemma}\label{Lemma: Characterizing morphisms}
Given two objects $\pair{\topo{G},U,N}$ and $\pair{\topo{G},V,M}$ in $\Eqsheav{G_1}{G_0}$, morphisms $\hat{t}:d^{-1}(U)/_{\sim_N}\rightarrow d^{-1}(V)/_{\sim_M}$ between them
%
%
%
%
are in one-to-one correspondence with open subsets
$T\subseteq d^{-1}(V)$
that satisfy the following properties:
\begin{enumerate}[i)]
\item $m(T\times_{G_0}M)\subseteq T$, i.e., $T$ is closed under $\sim_{M}$;
\item $c(T)=U$;
\item $m(T^{-1}\times_{G_0}T)\subseteq M$, i.e., if two arrows in $T$ share a codomain then they are $\sim_{M}$-equivalent;
\item $m(N\times_{G_0}T)\subseteq T$, i.e., if $f:x\rightarrow y$ is in $T$ and $g:y\rightarrow z$ is in $N$ then $g\circ f\in T$.
\end{enumerate}
Moreover, $\hat{t}$ can be thought of as `precomposing with $T$', in the sense that $\hat{t}([f]_{\sim_N})=[f\circ g]_{\sim_M}$ for some (any) $g\in T$ such that $c(g)=d(f)$.
\begin{proof}Straightforward.
\end{proof}
\end{lemma}
The following corollary will be useful.
\begin{corollary}\label{Corollary: Partial maps}
Given two objects of \cat{G}, \pair{\topo{G},U,N} and \pair{\topo{G},V,M}, and suppose $T\subseteq d^{-1}(V)$ is an open subset satisfying conditions (i), (iii), and (iv) of Lemma \ref{Lemma: Characterizing morphisms} and such that
$ ii')\ c(T)\subseteq U$
Then $T$ determines a morphism from the subobject \pair{\topo{G},c(T),N\upharpoonright_{c(T)}} of \pair{\topo{G},U,N} to \pair{\topo{G},V,M}.
\begin{proof}$c(T)$ is closed under $N$ by condition (iv), and the rest is straightforward.
\end{proof}
\end{corollary}
%
%
%
%
%
For a morphism $\fratil{f}{\topo{H}}{\topo{G}}$ of open topological groupoids,
%
%
%
%
%
%
%
%
%
the induced inverse image $f^*$ does not necessarily restrict to a functor between the respective Moerdijk-sites. The following condition (somewhat simplified from \cite{moerdijk:88}, cf. Lemma 6.2 there, so a proof is included here) ensures that it does.
\begin{definition}\label{Definition: Strictly full}
A morphism $\fratil{f}{\topo{H}}{\topo{G}}$ of open topological groupoids is a \emph{fibration} if for all $(h:x\rightarrow f_0(y))\in G_1$ there exists $g\in H_1$ such that $c(g)=y$ and $f_1(g)=h$. If the component continuous functions of $f$ are, moreover, subspace inclusions, then we say that \topo{H} is a \emph{replete subgroupoid} of \topo{G} and that $f$ is a \emph{replete subgroupoid inclusion}.
\end{definition}
Thus a replete subgroupoid is a full subcategory closed under isomorphisms and equipped with subspace topologies. Now, if $f:\topo{H}\rightarrow \topo{G}$ is a morphism of open topological groupoids and $(U,N)$ is an open subgroupoid of \topo{G}, then one readily sees that $(f_0^{-1}(U),f_1^{-1}(N))$ is an open subgroupoid of \topo{H}. Moreover:
%
\begin{lemma}\label{Lemma: Strictly full gives morphism of Moerdijk sites}
Let $\fratil{f}{\topo{H}}{\topo{G}}$ be a fibration of open topological groupoids, and let \pair{\topo{G},U,N} be an object of the Moerdijk-site of \Eqsheav{G_1}{G_0}. Then
\[\pair{\topo{H},f_0^{-1}(U),f_1^{-1}(N)}\cong f^*(\pair{\topo{G},U,N})\]
Moreover, if
\[\hat{t}:\pair{\topo{G},U_1,N_1}\rightarrow \pair{\topo{G},U_2,N_2}\]
is a morphism  in the Moerdijk-site of \Eqsheav{G_1}{G_0} corresponding to an open set $T\subseteq G_1$. Then
\[f^*(\hat{t}):\pair{\topo{H},f_0^{-1}(U_1),f_1^{-1}(N_1)}\rightarrow \pair{\topo{H},f_0^{-1}(U_2),f_1^{-1}(N_2)}\]
corresponds to the open set $f_1^{-1}(T)\subseteq H_1$.
\begin{proof}
Consider the diagram
\[\bfig
\square|alrm|/>`<-`>`>/<1200,400>[d^{-1}(f_0^{-1}(U))/_{\sim_{N_t}}`H_0\times_{G_0}d^{-1}(U)/_{\sim_N}`V=f_0^{-1}(U)`H_0;\hat{t}`e``\subseteq]
\place(1300,300,)[\pbangle]
\morphism<1200,400>[V=f_0^{-1}(U)`H_0\times_{G_0}d^{-1}(U)/_{\sim_N};t]
\square(1200,0)<1200,400>[H_0\times_{G_0}d^{-1}(U)/_{\sim_N}`d^{-1}(U)/_{\sim_N}`H_0`G_0;``c`f_0]
\efig\]
where $t$ is the section obtained by pulling back the canonical section $e:U\rightarrow d^{-1}(U)/_{\sim_N}$---so that $t(v)=\pair{v,[1_{f_0(v)}]_{\sim_N}}$---and $N_t\subseteq H_1$ and $\hat{t}$ are the induced open subgroupoid and morphism as in (\ref{Equation: Nt and t-hat}) and Proposition \ref{Proposition: GUNS generate}. Now, we have
\begin{align*}
N_t&=\cterm{g\in d^{-1}(V)\cap c^{-1}(V)}{f_1(g)\circ [1_{f_0(d(g))}]_{\sim_N}=[1_{f_0(c(g))}]_{\sim_N}}\\
   &=\cterm{g\in d^{-1}(V)\cap c^{-1}(V)}{f_1(g)\in N}=f_1^{-1}(N)
\end{align*}
and so $\pair{\topo{H},f_0^{-1}(U),f_1^{-1}(N)}=\pair{\topo{H},f_0^{-1}(U),N_t}$, and by Proposition \ref{Proposition: GUNS generate}, $\hat{t}$ is injective. Remains to show that it is also surjective. Let \pair{x,[g:u\rightarrow f_0(x)]_{\sim_N}} be given. Since $\fratil{f}{\topo{H}}{ \topo{G}}$ is a fibration, there exist $(h:y\rightarrow x)\in H_1$ such that $f_1(h)=g$, and since, accordingly, $f_0(y)=u$ we have $h\in d^{-1}(f_0^{-1}(U))$. But then
\[\hat{t}([h]_{\sim_{N_t}})=\pair{x,f_1(h)\circ[1_{f_0(y)}]_{\sim_N}}=\pair{x,[f_1(h)]_{\sim_N}}=\pair{x,[g]_{\sim_N}}.\]
The second claim is a similar computation using Lemma \ref{Lemma: Characterizing morphisms}.
\end{proof}
\end{lemma}
\subsection{Subgroupoids and subtoposes}
\label{Subsection: Subgroupoids and Subtoposes}

Let \topo{G} be an open topological groupoid and $\inklusjon{\iota}{\topo{H}}{\topo{G}}$ a replete subgroupoid, that is, \topo{H} is a topological groupoid consisting of subspaces $H_1\subseteq G_1$ and $H_0\subseteq G_0$ such that $H_0$ is closed under isomorphisms in \topo{G} and the inclusions form a morphism of groupoids which is full as a functor. It follows that \topo{H} is an open groupoid. By Lemma \ref{Lemma: Strictly full gives morphism of Moerdijk sites}, the induced inverse image functor $\fratil{\iota^*}{\Eqsheav{G_1}{G_0}}{ \Eqsheav{H_1}{H_0}}$ restricts to a functor between the respective Moerdijk sites $\fratil{I}{\cat{S}_{\topo{G}}}{\cat{S}_{\topo{H}}}$. It is shown in this section that this functor is essentially full and essentially surjective, whence the geometric morphism $f$ is an inclusion of toposes.
%

Say, for present purposes, that a functor $\fratil{F}{\cat{C}}{\cat{D}}$ is \emph{essentially full} if for any $B,C$ in \cat{C} and morphism $f:F(B)\rightarrow F(C)$ in \cat{D},
there exists in \cat{C} an object $B'$ with a zig-zag between $B$ and $B'$, and object $C'$ with a zig-zag between $C$ and $C'$, and a morphism $f':B'\rightarrow C'$ such that: i) $F$ sends the morphisms in both  zig-zags to  isomorphisms; and ii) the resulting isomorphisms $F(B)\cong F(B')$ and $F(C)\cong F(C')$ form a  commuting square with $f$ and $F(f)$:
\[\bfig \square<500,300>[F(B')`F(C')`F(B)`F(C)';F(f')`\cong`\cong`f] \efig\]
%
%
%
%
The following lemma will be useful.
\begin{lemma}\label{Lemma: Why we have strict fullness}
Let $\topo{H}\to/^{ (}->/<150>\topo{G}$ be a replete subgroupoid of an open groupoid, and let $V,W\subseteq G_1$ be open sets. Then $m(V\times_{G_0}W)$ is open and
\[m(V\times_{G_0}W)\cap H_1= m(V\cap H_1\times_{H_0}W\cap H_1)\]
\begin{proof}Composition of arrows is an open map for all open groupoids (see \cite{moerdijk:88}).
The rest is a straightforward consequence of the inclusion being  a fibration.
\end{proof}
\end{lemma}
\begin{lemma}\label{Lemma: I is essentially surjective}Let $\inklusjon{\iota}{\topo{H}}{\topo{G}}$ be a replete subgroupoid inclusion of open topological groupoids. Then the induced functor
$\fratil{I}{\cat{S}_{\topo{G}}}{\cat{S}_{\topo{H}}}$ of Moerdijk sites is essentially surjective and essentially full.
\begin{proof}
Consider an object \pair{\topo{H},V,M}. With $M$ an open set in the subspace $H_1\subseteq G_1$, we have the open sets
\[N:=\bigcup\cterm{K\in \cat{O}(G_1)}{K\cap H_1\subseteq M}\subseteq G_1\]
(where $\cat{O}(G_1)$ is the frame of open subsets of $G_1$) and $ U:=d(N)\cup c(N)$. Using Lemma \ref{Lemma: Why we have strict fullness}, it is straightforward to verify that $(U,N)$ is an open subgroupoid of \topo{G}, and clearly
%
%
%
%
%
%
%
$I(\pair{\topo{G},U,N})=\pair{\topo{H},U\cap H_0,N\cap H_1}=\pair{\topo{H},V,M}$. Thus this construction results in a right inverse $J$ to the the object function $I_0$.

Next, let \pair{\topo{G},U,N} be given. There is a canonical  morphism $\pair{\topo{G},U,N}\rightarrow J(I(\pair{\topo{G},U,N}))$ such that $I$ sends this morphism to the identity: Write
\begin{align*}
\pair{\topo{H},\underline{U},\underline{N}}&:=\pair{\topo{H},U\cap H_0,N\cap H_1}=I(\pair{\topo{G},U,N})\\
\pair{\topo{G},\overline{U},\overline{N}}&:=J(\pair{\topo{H},\underline{U},\underline{N}})
\end{align*}
Then $N\subseteq \overline{N}$ and, consequently, $U\subseteq \overline{U}$. Compose the canonical section $e:\overline{U}\rightarrow d^{-1}(\overline{U})$ with the inclusion $U\subseteq \overline{U}$,
\begin{equation}\label{Equation: Inproofref1}\bfig
\dtriangle|blb|/<-`>`>/<500,350>[d^{-1}(\overline{U})/_{\sim_{\overline{N}}}`\overline{U}`G_0;e`c`\subseteq]
\ptriangle(-500,0)|ala|/>`>`<-/<1000,350>[d^{-1}({U})/_{\sim_{{N}}}`d^{-1}(\overline{U})/_{\sim_{\overline{N}}}`U;\hat{v}`c`v]
\morphism(-500,0)|b|<500,0>[U`\overline{U};\subseteq]
%
\efig\end{equation}
For any $f:x\rightarrow y$ in $N$, we have that \[f\circ v(x)=f\circ e(x)=f\circ [1_x]_{\sim_{\overline{N}}}=[f]_{\sim_{\overline{N}}}=[1_y]_{\sim_{\overline{N}}}=v(y)\] since $N\subseteq\overline{N}$. So $v$ induces the morphism $\hat{v}([f]_{\sim_{N}})=[f]_{\sim_{\overline{N}}}$ in (\ref{Equation: Inproofref1}). By Lemma \ref{Lemma: Strictly full gives morphism of Moerdijk sites},
$\hat{v}$ is sent to the identity by $I$.
Now, given objects \pair{\topo{G},U,N}, \pair{\topo{G},V,M} and a morphism $\hat{t}:\pair{\topo{H},\underline{U},\underline{N}}\rightarrow\pair{\topo{H},\underline{V},\underline{M}}$, write $T\subseteq d^{-1}(\underline{V})$ for the corresponding open subset of arrows and $\hat{v}:\pair{\topo{G},{V},{M}} \rightarrow\pair{\topo{G},\overline{V},\overline{M}}$ for the morphism of the preceding paragraph.
Consider the open set
\[S:=c^{-1}(U)\cap\bigcup\cterm{P\in \cat{O}(G_1)}{P\cap H_1\subseteq T}.\]
It is straightforward to verify that $S$ satisfies the conditions of Corollary \ref{Corollary: Partial maps} so  that $S$ corresponds to a morphism $\hat{s}: \pair{\topo{G},c(S),N\upharpoonright_{c(S)}}\rightarrow\pair{\topo{G},\overline{V},\overline{M}}$,
\[\bfig \square/>`>`<-`/<1000,300>[\pair{\topo{G},c(S),N\upharpoonright_{c(S)}}`\pair{\topo{G},\overline{V},\overline{M}}`\pair{\topo{G},U,N}`\pair{\topo{G},V,M};\hat{s}`\subseteq`\hat{v}`]  \efig\]
where (by inspection and the proof of Lemma \ref{Lemma: I is essentially surjective}, respectively) $I$ sends both vertical arrows to identities. Moreover, $S\cap H_1=T$ and so by Lemma \ref{Lemma: Strictly full gives morphism of Moerdijk sites}, $I(\hat{s})=\hat{t}$.
\end{proof}
\end{lemma}
In conclusion:
\begin{theorem}\label{Theorem: Subgroupoids induce subtoposes}
Let \topo{G} be an open groupoid and $\fratil{\iota}{\topo{H}}{\topo{G}}$ a replete subgroupoid inclusion. Then the induced geometric morphism
\[\fratil{\iota}{\Eqsheav{H_1}{H_0}}{\Eqsheav{G_1}{G_0}}\]
is an inclusion.
\begin{proof}
By Lemma \ref{Lemma: I is essentially surjective} the inverse image $\fratil{\iota^*}{ \Eqsheav{G_1}{G_0}}{\Eqsheav{H_1}{H_0}}$ restricts to an essentially surjective and essentially full functor $\fratil{I}{\cat{S}_{\topo{G}}}{\cat{S}_{\topo{H}}}$. Consider the surjection-inclusion factorization
\[\bfig
\Vtriangle/>`->>`<-_{)}/<500,300>[\Sh{\topo{H}}`\Sh{\topo{G}}`\cat{I};\iota`e`m]
\efig\]
of $\iota$. The full subcategory $S_{\cat{I}}\to/^{ (}->/<150> \cat{I}$ consisting of the objects that are in  $m^*(\cat{S}_{\topo{G}})$ is a site for \cat{I} when equipped with the canonical coverage inherited from \cat{I}.
The inverse image $e^*$ restricts to a conservative functor $\fratil{E}{S_{\cat{I}}}{\cat{S}_{\topo{H}}}$ such that a family of morphisms in $S_{\cat{I}}$ is covering if and only if the image of it under $E$ is covering in $\cat{S}_{\topo{G}}$. But now $E$ is also essentially surjective, because $I$ is, and full, because it reflects isomorphisms and $F$ is essentially full. So $e$ is an equivalence.
\end{proof}
\end{theorem}
Note that in the special case where $H_0$ is an \emph{open} subset of $G_0$ (equivalently, $H_1$ is an open subset of $G_1$) the theorem follows from \cite[Prop.\ 5.13]{moerdijk:88} or from observing that in that case \topo{H} can be considered as a subterminal object of \Eqsheav{G_1}{G_0}. We shall return to this special case in Proposition \ref{Proposition: Open and closed inclusions} below.

\section{Quotient Theories and Subgroupoids}
\label{Section: Quotient Theories and Sub-Groupoids}
\subsection{Subtoposes, Quotient Theories, and Subgroupoids}
\label{Subsection: Quotient theories and subtoposes}
Let $\Sigma$ be a (first-order) signature. A geometric formula over $\Sigma$ is one constructed with the logical constants $\top$, $\bot$, $\wedge$, $\exists$, and $\bigvee$ (where the latter is infinitary disjunction of formulas that together have only finitely many free variables). See Part D of \cite{elephant1} for further details and a calculus for geometric sequents, i.e.\ sequents consisting of geometric formulas. Strictly speaking, there is a proper class of geometric formulas over $\Sigma$, but every geometric formula is provably equivalent, in the empty theory, to a disjunction of regular formulas (built from $\top$, $\wedge$, and $\exists$). We will therefore allow ourselves to speak of e.g.\ the collection \cat{L} of all sequents over $\Sigma$ as set instead of a class. It is convenient for our purposes to stipulate that theories are always closed under consequence, so
by a geometric theory is meant a deductively closed set of geometric sequents.
For theories \theory\ and $\theory'$ over the same signature $\Sigma$,  say that $\theory'$ is a \emph{quotient} of \theory\ and write
$\theory\subseteq\theory'$
if \theory\ is contained in $\theory'$ as a set of sequents.
Quotient theories of a theory \theory\ correspond to subtoposes of the classifying topos \classtop{T} (see \cite[Theorem 3.6]{caramello:99LoT}). Specifically, let \theory\ be a geometric theory. Recall from e.g.\ \cite{elephant1} that its classifying topos can be constructed by taking sheaves on the (essentially small) geometric syntactic category of \theory\ equipped with the coverage consisting of all sieves generated by small covering families
\[\classtop{T}:=\Sh{\synt{C}{T},J}\]
The subtoposes of \classtop{T} are then in 1--1 correspondence with the coverages on \synt{C}{T} containing $J$. Ordering coverages  by inclusion, this is an order reversing isomorphism of posets. Furthermore, the coverages containing $J$ are in 1--1 correspondence with the quotient theories of \theory, forming an isomorphism of posets when ordering quotient theories by inclusion. More details and further analysis regarding this correspondence can be found in \cite{caramello:99LoT}.

On the `geometric' side of things, consider an open topological groupoid \topo{G} and its equivariant sheaf topos $\Eqsheav{G_1}{G_0}$. If $\mathrm{Sub}(\Eqsheav{G_1}{G_0})$ is the poset of subtoposes of \Eqsheav{G_1}{G_0} and $\mathrm{Sub}(\topo{G})$ is the poset of replete subgroupoids of \topo{G} (isomorphic to the set of replete subsets of $G_0$ ordered by inclusion), then Theorem \ref{Theorem: Subgroupoids induce subtoposes} yields a morphism of posets
\[\fratil{\mathrm{sh}}{\mathrm{Sub}(\topo{G})}{\mathrm{Sub}(\Eqsheav{G_1}{G_0})}.\]
Now, suppose $\cat{F}$ is a subtopos of \Eqsheav{G_1}{G_0}. Since every element $x$ of $G_0$ induces a point   $\fratil{p_x}{\Sets}{\Eqsheav{G_1}{G_0}}$ (and every element of $G_1$ an invertible geometric transformation of points), we can form the (replete) subset $H_0\subseteq G_0$ of those elements that induce points that factor through \cat{F}. This yields a morphism
\[\fratil{\mathrm{pt}}{\mathrm{Sub}(\Eqsheav{G_1}{G_0})}{\mathrm{Sub}(\topo{G})}.\]
There is, accordingly, a connection between quotients of the theory classified by \Eqsheav{G_1}{G_0} and subgroupoids of \topo{G}, which we state next together with a characterization of the subgroupoids in the image of \mng{pt}.
For more on the general method of using the various ways in which toposes can be viewed and presented to mediate between different structures and theories see \cite{caramello:00TUMVTT}.

%
%

%
\subsection{Groupoids of Models and Definable Subsets}
\label{Subsection: Groupoids of Models}
Let \topo{G} be a topological groupoid, \Eqsheav{G_1}{G_0} the topos of equivariant sheaves on it. Then (see \cite{elephant1}) there exists a geometric theory, \theory, such that
\begin{equation}\label{Equation: groupoid-theory}\Eqsheav{G_1}{G_0}\simeq \classtop{T}\simeq \Sh{\synt{C}{T},J}\end{equation}
Since an element of $G_0$ induces a point of this topos, and an element of $G_1$ induces an invertible geometric transformation of points, we can, by the equivalence between the category of \theory-models and the category of points of \classtop{T}, regard \topo{G} as a topological groupoid of \theory-models and isomorphisms. $G_0$ is then a space of \emph{enough} models for \theory, in the sense that is a sequent is true in all models in $G_0$ then it is in \theory. This follows since the points induced by elements of $G_0$ are \emph{enough} for \Eqsheav{G_1}{G_0} in the sense that the inverse image functors of the induced points are jointly conservative (see \cite{elephant1}). Conversely, given a theory \theory\ with enough models, \cite{butz:98b} constructs an open topological groupoid of models and isomorphisms such that  $\Eqsheav{G_1}{G_0}\simeq \classtop{T}$. (More direct---in logical terms---variations of this construction for geometric and classical first-order theories respectively can also be found in \cite{MathQuart:repforgeothrs} and the  \cite{APAL:fol}.)

Fix an open topological groupoid \topo{G}, a theory \theory\ over a signature $\Sigma$, and an equivalence as displayed in (\ref{Equation: groupoid-theory}) above, and regard \topo{G} as a groupoid of \theory-models and isomorphisms accordingly. We shall write elements as $\alg{M}, \alg{N}\in G_0$ and $\alg{f}, \alg{g}\in G_1$ when we want to emphasize this perspective.
%
\begin{lemma}\label{Lemma: Galois translated}
\textup{(}i\textup{)} Let $\theory'$ be a quotient theory of \theory\ and \classtop{T'} the corresponding subtopos of \Eqsheav{G_1}{G_0}. Then $\mng{pt}(\classtop{T'})=\topo{H}$ where
\[H_0=\cterm{\alg{M}\in G_0}{\alg{M}\vDash\theory'}\subseteq G_0.\]
\textup{(}ii\textup{)} Let \topo{H} be a subgroupoid of  \topo{G}. Then $\mng{sh}(\topo{H})$ classifies the quotient theory
\[\theory'=\cterm{\sigma\in\cat{L}}{\alg{M}\vDash\sigma,\ \textnormal{for all}\ \alg{M}\in H_0}\supseteq \theory\]
where \cat{L} is the set of all geometric sequents over $\Sigma$.
\begin{proof}(i) A point \fratil{p_{\alg{M}}}{\Sets}{\Eqsheav{G_1}{G_0}} induced by $\alg{M}\in G_0$ factors through the subtopos \classtop{T'} if and only if $\alg{M}\vDash \theory'$.

\noindent (ii) Let $\theory'$ be the quotient theory classified by $\mng{sh}(\topo{H})$. Clearly, $\alg{M}\vDash\theory'$ for all $\alg{M}\in H_0$. Since the points induced by elements of $H_0$ are enough for  \Eqsheav{H_1}{H_0}, it is also the case that if $\sigma$ is a sequent true in all models in $H_0$, then $\sigma\in\theory'$. Thus the quotient $\theory'$ is determined by the subset $H_0$ as the set of sequents true in all models in $H_0$.
\end{proof}
\end{lemma}
\begin{proposition}\label{proposition: Galois connection}
Let \topo{G} be an open topological groupoid. The morphisms of posets $\mng{pt}: \mathrm{Sub}(\Eqsheav{G_1}{G_0})\leftrightarrows \mathrm{Sub}(\topo{G}) :\mng{sh}$ form a Galois connection
\[\begin{array}{ccc}
\mng{sh}(\topo{H})&\leq &\cat{F}\\[3pt]
\hline\\[-2.6ex]
\hline
\\[-9pt]
\topo{H}&\leq &\mng{pt}(\cat{F})
\end{array}\]
between subtoposes of \Eqsheav{G_1}{G_0} and subgroupoids of \topo{G}.
\begin{proof}By the subtopos-quotient theory correspondence, since it is clear from Lemma \ref{Lemma: Galois translated} that the  quotient theory classified by $\cat{F}\hookrightarrow\Eqsheav{G_1}{G_0}$ is contained in the quotient theory classified by $\mng{sh}(\mng{pt}(\cat{F}))$.
\end{proof}
\end{proposition}
Say that an open topological groupoid \topo{G} is \emph{saturated} (with apologies for overloading that term) if every subtopos of \Eqsheav{G_1}{G_0} with enough points is of the form \Eqsheav{H_1}{H_0} for a subgroupoid $\topo{H}\to/^{ (}->/<150> \topo{G}$; equivalently, if every subtopos with enough points has enough points induced by elements of $G_0$. In logical terms, with respect to a classified theory \theory\ as in (\ref{Equation: groupoid-theory}), this is saying that for any quotient theory $\theory'$ of \theory, if $\theory'$ has enough models, the models in the set $G_0$ are already enough. Since the groupoids of models and isomorphisms constructed in \cite{butz:98b} (and \cite{MathQuart:repforgeothrs} and \cite{APAL:fol}) are by their construction saturated in this sense, we restrict attention to saturated groupoids. Say that a subgroupoid is \emph{definable} if it is in the image of \mng{pt}, or from a logical perspective, if it is of the form $\cterm{\alg{M}\in G_0}{\alg{M}\vDash\theory'}\subseteq G_0$ for a quotient theory $\theory'$ of \theory.
We proceed to characterize the definable subgroupoids of a saturated groupoid \topo{G} directly in terms of the groupoid.

\begin{definition}\label{Definition: Geometric domination}
For an open topological groupoid \topo{G}, an element $x\in G_0$ and a subset $H_0\subseteq G_0$, say that $x$ (geometrically) \emph{dominates} $H_0$, written $x\gg_{\small GD}H_0$, if for all open subgroupoids (U,N) of \topo{G} and all open subsets $V,W\subseteq U$ that are closed under $N$ we have
\begin{align*}
c^{-1}(H_0)\cap d^{-1}(V)\subseteq  d^{-1}(W)\\
\Rightarrow\ c^{-1}(x)\cap d^{-1}(V)\subseteq  d^{-1}(W)
\end{align*}
\end{definition}
\begin{theorem}\label{Theorem: Definable subgroupoids}
Let \topo{G} be an open
topological groupoid, and \topo{H} a replete subgroupoid. Then \topo{H} is definable iff $H_0$ is closed under domination, in the sense that for any $x\in G_0$ if $x\gg_{GD}H_0$ then $x\in H_0$.
\begin{proof}
Let \theory\ be a geometric theory such that $\Eqsheav{G_1}{G_0}\simeq\classtop{T}$. Then, corresponding to to the generic model,  we can choose a small, geometric, full subcategory $\cat{T}$ of $\Eqsheav{G_1}{G_0}$ (closed under subobjects) the objects of which form a generating set.  On the other hand, $\cat{S}_{\topo{G}}$ is  a small, full subcategory (closed under subobjects) the objects of which form a generating set. Write $p_x$ for the point induced by $x\in G_0$. Then saying that $H_0$ is definable comes to saying that for all $x\in G_0$, if for all objects $A\in \cat{T}$ and subobjects $P,Q\rightarrowtail A$, if $p_x^*(P)\leq p_x^*(Q)$ whenever $p_y^*(P)\leq p_y^*(Q)$ for all $y\in H_0$, then $x\in H_0$. Similarly, by Lemma \ref{Lemma: Subobjects in the Moerdijk site}, saying that $H_0$ is closed under domination comes to saying that for all $x\in G_0$, if for all objects $A\in \cat{S}_{\topo{G}}$ and subobjects $P,Q\rightarrowtail A$, if $p_x^*(P)\leq p_x^*(Q)$ whenever $p_y^*(P)\leq p_y^*(Q)$ for all $y\in H_0$, then $x\in H_0$. Since \cat{T} and $\cat{S}_{\topo{G}}$ are generating, this is equivalent.

\end{proof}
\end{theorem}
%


 Lemma \ref{Lemma: Galois translated}, Proposition \ref{proposition: Galois connection}, and Theorem \ref{Theorem: Definable subgroupoids} open up the possibility of extending the analysis of the correspondence between quotient theories and subtoposes to include subgroupoids. For instance, \cite{caramello:99LoT} contains detailed proofs that open subtoposes correspond to quotient theories obtained by adding a single geometric sentence as an axiom, and closed subtoposes to quotient theories obtained by adding a single sequent of the form $\phi\vdash\bot$ where $\phi$ is a geometric sentence. In terms of subgroupoids, we have the following.
\begin{proposition}\label{Proposition: Open and closed inclusions}
Let $\topo{H}$ be a definable subgroupoid of an open, saturated topological groupoid \topo{G}, and fix \theory\ such that  \Eqsheav{G_1}{G_0} classifies \theory.
\begin{enumerate}
\item \Eqsheav{H_1}{H_0} classifies a quotient $\theory'$ such that $\theory'$ can be obtained from \theory\ by adding a single geometric sentence as an axiom if and only if $H_0\subseteq G_0$ is an open subset.
\item \Eqsheav{H_1}{H_0} classifies a quotient $\theory'$ such that $\theory'$ can be obtained from \theory\ by adding a single geometric sequent $\phi\vdash\bot$ as an axiom where $\phi$ is a geometric sentence if and only if  $H_0\subseteq G_0$ is a closed subset.
\end{enumerate}
\begin{proof}

(1) As noted, \Eqsheav{H_1}{H_0} classifies a quotient $\theory'$ such that $\theory'$ can be obtained from \theory\ by adding a single geometric sentence if and only if \Eqsheav{H_1}{H_0} is an open subtopos. If $H_0\subseteq G_0$ is open and closed under $G_1$, we can consider $H_0$ as a subterminal object, slicing over which produces the (inverse image part of) the induced geometric inclusion, which is thereby open. Conversely,  the (inverse image part of) the induced geometric inclusion is up to equivalence obtained by slicing over a subterminal object, and a subterminal object can be considered as an open subset $U\subseteq G_0$ closed under $G_1$. Now, $U$ must be definable---i.e.\ closed under domination---for if $x\gg_{\textnormal{\tiny GD}}U$ then
\begin{align*}
c^{-1}(U)\cap d^{-1}(G_0)\subseteq  d^{-1}(U)\\
\Rightarrow\ c^{-1}(x)\cap d^{-1}(G_0)\subseteq  d^{-1}(U)
\end{align*}
implies that $x\in U$. But then $U=H_0$ since both are definable and they classify the same theory.

(2) By the above, $H_0\subseteq G_0$ is closed if and only if  there exists a single geometric sentence $\phi$ such that $H_0$ is the set of $\theory$-models (in $G_0$) where $\phi$ is false if and only if  $H_0$ is defined by the theory (generated by) $\theory\cup\{\phi\vdash\bot\}$ for a geometric sentence $\phi$ (note that if a theory has enough models, then so does any quotient obtained by adding a single sequence of the form $\phi\vdash\bot$ for a sentence $\phi$).
\end{proof}
\end{proposition}
%

\bibliographystyle{ieeetr}
\bibliography{bibliografi}

\begin{thebibliography}{1}

\bibitem{butz:98b}
C.~Butz and I.~Moerdijk, ``Representing topoi by topological groupoids,'' {\em
  Journal of Pure and Applied Algebra}, vol.~130, pp.~223--235, 1998.

\bibitem{caramello:99LoT}
O.~Caramello, ``Lattices of theories,'' 2009.
\newblock {\tt http://arxiv.org/PS\_cache/arxiv/pdf/0905/0905.0299v1.pdf}.

\bibitem{moerdijk:88}
I.~Moerdijk, ``The classifying topos of a continuous groupoid. {I},'' {\em
  Transactions of the American Mathematical Society}, vol.~310, no.~2,
  pp.~629--668, 1988.

\bibitem{elephant1}
P.~T. Johnstone, {\em Sketches of an Elephant}, vol.~43 and 44 of {\em Oxford
  Logic Guides}.
\newblock Oxford: Clarendon Press, 2002.

\bibitem{moerdijk:90}
I.~Moerdijk, ``The classifying topos of a continuous groupoid. {I}{I},'' {\em
  Cahiers de Topologie et G\'eom\'etrie Diff\'erentielle Cat\'egoriques},
  vol.~XXXI, no.~2, pp.~137--167, 1990.

\bibitem{onlinenote:subgpds}
H.~Forssell, ``Subgroupoids and quotient theories.''
\newblock {\tt http://arxiv.org/abs/1111.2952}.

\bibitem{caramello:00TUMVTT}
O.~Caramello, ``The unification of mathematics via topos theory,'' 2010.
\newblock {\tt http://arxiv.org/abs/1006.3930}.

\bibitem{MathQuart:repforgeothrs}
H.~Forssell, ``Topological representation of geometric theories,'' {\em
  Mathematical Logic Quarterly}, vol.~58, pp.~380--393, 2012.

\bibitem{APAL:fol}
S.~Awodey and H.~Forssell, ``First-order logical duality,'' {\em Annals of Pure
  and Applied Logic}, vol.~164, pp.~319--348, 2013.

\end{thebibliography}
\end{document}